\documentclass[a4paper,12pt]{article}
\usepackage[utf8]{inputenc}
\usepackage[T1]{fontenc}
\usepackage{amsmath, amssymb, amsthm}
\usepackage{graphicx}
\usepackage{hyperref}
\usepackage[margin=1in]{geometry}
\usepackage{listings}
\usepackage{xcolor}
\usepackage{enumitem}
\usepackage{fancyhdr}
\usepackage{booktabs}
\usepackage{caption}
\usepackage{subcaption}

\newcommand{\squaternion}{\mathbf{q}}
\newcommand{\splitq}{\mathbb{S}}
\newcommand{\kay}{k}
\newcommand{\conj}[1]{\overline{#1}}

\lstdefinestyle{matlab}{
    language=Matlab,
    basicstyle=\ttfamily\small,
    keywordstyle=\color{blue},
    commentstyle=\color{green!60!black},
    stringstyle=\color{red},
    numbers=left,
    numberstyle=\tiny\color{gray},
    stepnumber=1,
    numbersep=10pt,
    backgroundcolor=\color{white},
    showspaces=false,
    showstringspaces=false,
    showtabs=false,
    frame=single,
    tabsize=4,
    captionpos=b,
    breaklines=true,
    breakatwhitespace=false,
    title=\lstname,
    escapeinside={(*}{*)},
    morekeywords={SplitQuaternionMatrix, SplitQuaternion}
}

\hypersetup{
    colorlinks=true,
    linkcolor=blue,
    urlcolor=blue,
    citecolor=blue
}

\title{Split Quaternion Toolbox: A MATLAB Implementation 
}

 \author{Salman Ahmadi-Asl\\ \texttt{s.ahmadiasl@innopolis.ru}}
 \author{
     Salman Ahmadi-Asl \\
     \small Innopolis University, Innopolis, Russia \\
     \small \texttt{s.ahmadiasl@innopolis.ru}
      \and
      Abbas Moradi \\ 
      \small Statistical Research and Training Center (SRTC), Tehran, Iran \\
      \small \texttt{amoradi@srtc.ac.ir}
     \and
     Yang Zhang \\
     \small University of Manitoba, Winnipeg, Canada \\
     \small \texttt{yang.zhang@umanitoba.ca}
 }
 
 \date{}

\begin{document}

\maketitle

\begin{abstract}
This paper presents a MATLAB toolbox for working with split quaternions and split quaternion matrices. Split quaternions, also known as coquaternions, form a four-dimensional associative non-commutative algebra over the real numbers with a signature (+,+,-,-). The toolbox provides two main classes: SplitQuaternion for scalar-valued split quaternions and SplitQuaternionMatrix for matrix-valued split quaternions. The implementation includes arithmetic operations, algebraic functions, matrix operations, linear algebra utilities, and various conversion methods. The toolbox is designed to be intuitive for users familiar with MATLAB's standard matrix operations while providing the full power of split quaternion algebra. This paper details the implementation, provides numerous examples, and demonstrates the toolbox's capabilities through practical demonstrations. The MATLAB toolbox is accessible at \url{https://github.com/SalmanAhmadi-Asl/Split_Quaternion.git}.

\end{abstract}



\section{Introduction}

The study of split quaternion algebras has received a lot of attention in recent years because of their deep connections to theoretical physics, geometry, and numerical linear algebra. Unlike standard Hamilton quaternions, split quaternions have a hyperbolic signature that makes them essential for modeling phenomena in Minkowski space-time and Lorentzian geometry. Researchers have thoroughly explored the basic algebraic properties of split quaternions, including their roots and geometric interpretations. For example, \cite{ozdemir2009roots} gave a thorough analysis of the roots of a split quaternion, while \cite{ata2009split} looked at their role in semi-Euclidean projective spaces. \cite{atasoy2017new} moved the field forward by introducing a new polar representation for split and dual split quaternions, which adds to their geometric usefulness. More recently, \cite{ozturk2023geometric} have deepened our understanding of their geometric importance in applied settings.


Parallel to these algebraic advancements, the theory of split quaternion matrices has become an important research area. \cite{alagoz2012split} did early work on the basic properties of these matrices. Later, \cite{erdougdu2013complex} studied complex representations, and \cite{erdougdu2013eigenvalues} looked at the eigenvalue structures of split quaternion matrices, showing details that make them different from quaternion matrices. Researchers have also studied computational aspects; for example, \cite{jiang2018algebraic} developed algebraic methods for eigenvalue and eigenvector problems, which are important for applications in split quaternionic mechanics. \cite{wang2024singular} provided a strong theory of singular value decomposition (SVD) for split quaternion matrices, giving a tool that can be widely used in engineering and physics.

Matrix equations over split quaternion algebras form a core component of current research, driven by applications in control theory, signal processing, and quantum mechanics. The consistency and solvability of various classes of these equations have been investigated by \cite{liu2019consistency}, who studied equations of the form \( AX^\star - XB = CY + D \). \cite{yuan2017hermitian} focused on Hermitian solutions to the equation \( AXB + CXD = E \), while \cite{si2024classical} provided a comprehensive analysis of classical systems of matrix equations over the split quaternion algebra. The least squares problem, a cornerstone of numerical analysis, has also been adapted to this setting; \cite{liu2020least} explored least squares \( X = \pm X ^{\eta^*} \) solutions, and more recently, \cite{zhang2026efficient} developed efficient real algorithms based on the semi-tensor product (STP) for least squares problems in the complex split quaternion algebra.

The practical relevance of split quaternions is underscored by their applications in modern physics and signal processing. \cite{gogberashvili2014split} utilized split quaternions to describe particles in (2+1)-space, while \cite{guo2023novel} proposed a novel algebraic approach for solving the Schr\"{o}dinger equation within a split quaternionic mechanics framework. In a more applied domain, \cite{davila2023quaternion} demonstrated the efficacy of quaternion and split quaternion neural networks for low-light color image enhancement, highlighting the potential of this algebra in modern computational intelligence. Furthermore, \cite{ozdemir2021kinematic} explored the kinematic equations of Lorentzian magnetic flux tubes, bridging the gap between abstract algebra and differential geometry. The recent work by \cite{jia2025general} on the general solution to systems of tensor equations further extends the algebraic toolkit, promising new avenues for research. Iterative methods to compute the Moore-Penrose inverse of split quaternion matrices are studied in \cite{ahmadi2026iterative}. This growing body of work establishes split quaternion matrices and their associated equations as a rich and fertile ground for both theoretical exploration and practical application. 

There are software tools for regular quaternions (like \texttt{QTFM} and \texttt{Quaternion Toolbox for MATLAB}), but there is no complete tool for split quaternions. Our \texttt{SplitQuaternionToolbox} fills this gap by offering:
\begin{itemize}
    \item Basic split quaternion math operations
    \item Split quaternion matrix operations
    \item Advanced matrix methods
    \item Tools to convert between different formats
\end{itemize}

This paper is organized as follows: Section \ref{sec:pre} covers the basic definitions, properties, and matrix forms of split quaternions, and shows how they differ from regular quaternions. Section \ref{sec:toolbox_over} provides an overview of the toolbox. Several examples are presented in Section \ref{sec:exm} explaining how to use the toolbox. Section \ref{Sec:matrix_rep} deals with the matrix representation of split quaternions and split quaternion matrices. The CUR approximation of split quaternion matrices is discussed in Section \ref{SEc:CUR}. The conclusion is given in Section \ref{Sec:con}.

\section*{2. Preliminaries}\label{sec:pre}

This section establishes the algebraic and matrix-theoretic foundations of quaternions necessary for the developments that follow. We present both the scalar algebra and its matrix extension, emphasizing aspects relevant to computational linear algebra.


\subsection{Quaternions and quaternion matrices}

The quaternion algebra $\mathbb{H}$ extends complex numbers to four dimensions. An element $q \in \mathbb{H}$ is expressed as
\begin{equation}
    q = q_0 + q_1\mathbf{i} + q_2\mathbf{j} + q_3\mathbf{k},
    \label{eq:quaternion}
\end{equation}
with real coefficients $q_0, q_1, q_2, q_3 \in \mathbb{R}$. The imaginary units $\mathbf{i}, \mathbf{j}, \mathbf{k}$ satisfy the fundamental relations
\begin{equation}
    \mathbf{i}^2 = \mathbf{j}^2 = \mathbf{k}^2 = \mathbf{i}\mathbf{j}\mathbf{k} = -1,
    \label{eq:fundamental_relations}
\end{equation}
which imply the anticommutative pairwise products:
\begin{equation}
    \mathbf{i}\mathbf{j} = \mathbf{k} = -\mathbf{j}\mathbf{i}, \quad
    \mathbf{j}\mathbf{k} = \mathbf{i} = -\mathbf{k}\mathbf{j}, \quad
    \mathbf{k}\mathbf{i} = \mathbf{j} = -\mathbf{i}\mathbf{k}.
    \label{eq:anticommutative}
\end{equation}
This non-commutative multiplication defines the essential character of the algebra. The quaternion can be decomposed into its scalar part, $\operatorname{S}(q) = q_0$, and its vector part, $\operatorname{V}(q) = q_1\mathbf{i} + q_2\mathbf{j} + q_3\mathbf{k}$. Thus, $q = \operatorname{S}(q) + \operatorname{V}(q)$.


The basic algebraic operations are defined as follows.

\textbf{Conjugation.} The conjugate of $q$ is $\conj{q} = \operatorname{S}(q) - \operatorname{V}(q) = q_0 - q_1\mathbf{i} - q_2\mathbf{j} - q_3\mathbf{k}$. Conjugation is an anti-automorphism: $\conj{pq} = \conj{q} \,\conj{p}$ for any $p, q \in \mathbb{H}$.

\textbf{Norm.} The norm of $q$ is the non-negative real number $\|q\| = \sqrt{qq^*} = \sqrt{q_0^2 + q_1^2 + q_2^2 + q_3^2}$. A quaternion with unit norm is termed a \emph{versor}.

\textbf{Multiplication.} The product of two quaternions $p$ and $q$ can be compactly expressed using their scalar-vector forms:
\begin{equation}
    pq = \bigl(\operatorname{S}(p)\operatorname{S}(q) - \operatorname{V}(p) \cdot \operatorname{V}(q)\bigr) + \bigl(\operatorname{S}(p)\operatorname{V}(q) + \operatorname{S}(q)\operatorname{V}(p) + \operatorname{V}(p) \times \operatorname{V}(q)\bigr),
    \label{eq:product_formula}
\end{equation}
where $\cdot$ and $\times$ denote the standard Euclidean dot and cross products. This formula highlights the interplay between scalar multiplication and three-dimensional vector operations. To interface with standard numerical linear algebra, quaternions are often represented as real or complex matrices. These isomorphisms preserve the algebraic structure but expand the dimension. Below, we present the real and complex embeddings of quaternions:

\textbf{Real (${4\times4}$) Representation.}
Each quaternion $q$ corresponds to a linear transformation on $\mathbb{R}^4$ via left multiplication. This yields the real matrix representation:
\begin{equation}
    \mathcal{R}(q) =
    \begin{bmatrix}
         q_0 & -q_1 & -q_2 & -q_3 \\
         q_1 &  q_0 & -q_3 &  q_2 \\
         q_2 &  q_3 &  q_0 & -q_1 \\
         q_3 & -q_2 &  q_1 &  q_0
    \end{bmatrix}.
    \label{eq:real_rep}
\end{equation}
The map $\mathcal{R}: \mathbb{H} \to \mathbb{R}^{4\times4}$ is an injective homomorphism: $\mathcal{R}(pq) = \mathcal{R}(p)\mathcal{R}(q)$ and $\mathcal{R}(\conj{q}) = \mathcal{R}(q)^\top$.

\textbf{Complex (${2\times2}$) Representation.}
Exploiting the isomorphism $\mathbb{H} \cong \mathbb{C}^{2\times2}$, we have the complex Pauli representation:
\begin{equation}
    \mathcal{C}(q) =
    \begin{bmatrix}
        q_0 + iq_1 & q_2 + iq_3 \\
        -q_2 + iq_3 & q_0 - iq_1
    \end{bmatrix},
    \label{eq:complex_rep}
\end{equation}
where $i$ is the classical imaginary unit. This representation is particularly useful for eigenvalue analysis.


A quaternion matrix $\mathbf{A} = (a_{ij}) \in \mathbb{H}^{m\times n}$ has quaternion entries. The set $\mathbb{H}^{m\times n}$ forms a right module over $\mathbb{H}$, meaning scalar multiplication is defined from the right: $(\mathbf{A}\lambda)_{ij} = a_{ij}\lambda$. Operations are defined entry-wise, respecting the non-commutative base ring.

\textbf{Addition and Multiplication.} For $\mathbf{A}, \mathbf{B} \in \mathbb{H}^{m\times n}$, addition is component-wise. For $\mathbf{A} \in \mathbb{H}^{m\times p}$ and $\mathbf{B} \in \mathbb{H}^{p\times n}$, the product $\mathbf{C} = \mathbf{A}\mathbf{B}$ is defined by $c_{ij} = \sum_{k=1}^p a_{ik}b_{kj}$, maintaining the order of factors.

\textbf{Conjugate Transpose.} The conjugate transpose $\mathbf{A}^* \in \mathbb{H}^{n\times m}$ is given by $(\mathbf{A}^*)_{ij} = \conj{a_{ji}}$. This operation is an involution: $(\mathbf{A}^*)^* = \mathbf{A}$ and $(\mathbf{A}\mathbf{B})^* = \mathbf{B}^*\mathbf{A}^*$.

\textbf{Matrix Norms.} The Frobenius norm is $\|\mathbf{A}\|_F = \sqrt{\sum_{i,j} \|a_{ij}\|^2}$. This norm is submultiplicative and unitary invariant.


Several important classes of matrices emerge, mirroring their real and complex analogues but with quaternion-specific properties.

\begin{itemize}
    \item \textbf{Hermitian Matrices:} $\mathbf{A}^* = \mathbf{A}$. Their eigenvalues (in the right sense) are real numbers.
    \item \textbf{Skew-Hermitian Matrices:} $\mathbf{A}^* = -\mathbf{A}$.
    \item \textbf{Unitary Matrices:} $\mathbf{A}^*\mathbf{A} = \mathbf{A}\mathbf{A}^* = \mathbf{I}$. They preserve the quaternion norm: $\|\mathbf{A}\mathbf{x}\| = \|\mathbf{x}\|$ for $\mathbf{x} \in \mathbb{H}^n$.
\end{itemize}

Note that the real and complex representations of a quaternion matrix $\mathbf{A}=\mathbf{A}_1+\mathbf{A}_2\mathbf{i}+\mathbf{A}_3\mathbf{j}+\mathbf{A}_4\mathbf{k}$ will follow the same pattern as those in \eqref{eq:real_rep} and \eqref{eq:complex_rep}, where instead of $q_i$, we replace $\mathbf{A}_i$.

\subsection{Split quaternions and split quaternion matrices}

A split-quaternion is an element of the form:
\[
q = a + b\mathbf{i} + c\mathbf{j} + d\mathbf{k}\in\mathbb{S},
\]
where $a, b, c, d \in \mathbb{R}$ and the basis elements $\{\mathbf{1}, \mathbf{i}, \mathbf{j}, \mathbf{k}\}$ satisfy the following multiplication rules:
\begin{align*}
    \mathbf{i}^2 &= -1, \quad \mathbf{j}^2 = +1, \quad \mathbf{k}^2 = +1, \\
    \mathbf{ij} &= \mathbf{k} = -\mathbf{ji} \quad \text{(anti-commutativity)}, \\
    \mathbf{jk} &= -\mathbf{i} = -\mathbf{kj}, \\
    \mathbf{ki} &= \mathbf{j} = -\mathbf{ik}.
\end{align*}
which can be summarized in the following table:
\[
\begin{array}{c|cccc}
\times & \mathbf{1} & \mathbf{i} & \mathbf{j} & \mathbf{k} \\ \hline
\mathbf{1} & \mathbf{1} & \mathbf{i} & \mathbf{j} & \mathbf{k} \\
\mathbf{i} & \mathbf{i} & -\mathbf{1} & \mathbf{k} & -\mathbf{j} \\
\mathbf{j} & \mathbf{j} & -\mathbf{k} & +\mathbf{1} & -\mathbf{i} \\
\mathbf{k} & \mathbf{k} & \mathbf{j} & \mathbf{i} & +\mathbf{1} \\
\end{array}
\]

Note the mixed signature: one imaginary unit squares to $-1$ (as in complex numbers), while two square to $+1$ (as in split-complex or hyperbolic numbers). For $\squaternion = q_0 + q_1{\mathbf i} + q_2{\mathbf j} + q_3 {\mathbf k}$ and $\mathbf{p} = p_0 + p_1{\mathbf i} + p_2{\mathbf j} + p_3{\mathbf k}$, the basic operations are defined as follows:

\begin{itemize}
    \item \textbf{Addition:}
    \[
    \squaternion + \mathbf{p} = (q_0 + p_0) + (q_1 + p_1){\mathbf i} + (q_2 + p_2){\mathbf j} + (q_3 + p_3){\mathbf k}.
    \]
    
    \item \textbf{Multiplication:}
    \begin{align*}
    \squaternion\mathbf{p} =& (q_0p_0 - q_1p_1 + q_2p_2 + q_3p_3), \\
    &+ (q_0p_1 + q_1p_0 + q_2p_3 - q_3p_2){\mathbf i}, \\
    &+ (q_0p_2 - q_1p_3 + q_2p_0 + q_3p_1){\mathbf j}, \\
    &+ (q_0p_3 + q_1p_2 - q_2p_1 + q_3p_0){\mathbf k}.
    \end{align*}
    
    \item \textbf{Conjugate:}
    \[
    \conj{\squaternion} = q_0 - q_1{\mathbf i} - q_2{\mathbf j} - q_3{\mathbf k}.
    \]
    
    \item \textbf{Norm:} The split quaternion norm is given by:
    \[
    N(\squaternion) = \squaternion\conj{\squaternion} = \conj{\squaternion}\squaternion = q_0^2 + q_1^2 - q_2^2 - q_3^2.
    \]
    Note that this norm is not positive-definite but takes values in $\mathbb{R}$. Unlike ordinary quaternions, the norm $N(q)$ can be zero for non-zero $q$, making the algebra more challenging numerically. Based on the norm $N(\squaternion) = q_0^2 + q_1^2 - q_2^2 - q_3^2$, the following classifications are defined:
\begin{itemize}
    \item \textbf{Elliptic:} $N(\squaternion) > 0$.
    \item \textbf{Hyperbolic:} $N(\squaternion) < 0$.
    \item \textbf{Parabolic:} $N(\squaternion) = 0$, $\squaternion \neq 0$.
    \item \textbf{Zero:}$\squaternion = 0$.
\end{itemize}

    \item \textbf{Inverse:} If $N(\squaternion) \neq 0$, the inverse exists:
    \[
    \squaternion^{-1} = \frac{\conj{\squaternion}}{N(\squaternion)}.
    \]
\end{itemize}

The split quaternion algebra is isomorphic to the real $4 \times 4$ matrix algebra:
\begin{equation}\label{eq:real_split}
\squaternion = q_0 + q_1\imath + q_2\jmath + q_3\kay \cong 
\begin{pmatrix}
q_0  & -q_1 & q_2& q_3\\
q_1 & q_0 & q_3& -q_2\\
q_2&q_3& q_0& -q_1\\
q_3&-q_2&q_1&q_0
\end{pmatrix},
\end{equation}
and also the $2\times 2$ complex matrix algebra
\begin{equation}\label{eq:complex_split}
\squaternion = q_0 + q_1\imath + q_2\jmath + q_3\kay \cong 
\begin{pmatrix}
q_0 + iq_2 & q_1 + iq_3 \\
-iq_1 + q_3 & q_0 - iq_2
\end{pmatrix}.
\end{equation}
It is also not difficult to show that
\begin{equation}\label{eq:split_real_2}
\squaternion = q_0 + q_1\imath + q_2\jmath + q_3\kay \cong 
\begin{pmatrix}
q_0 + q_2 & q_1 + q_3 \\
-q_1 + q_3 & q_0 - q_2
\end{pmatrix}.
\end{equation}
Key properties and comparisons with quaternions are summarized in Table \ref{tab:quat_split_compare}. 

\begin{table}[htbp]
\centering
\small 
\caption{Comparison of algebraic properties: quaternions vs. split-quaternions}
\label{tab:quat_split_compare}
\begin{tabular}{p{0.2\textwidth}p{0.35\textwidth}p{0.35\textwidth}} 
\toprule
\textbf{Feature} & \textbf{Quaternions ($\mathbb{H}$)} & \textbf{Split-Quaternions} \\
\midrule
\textbf{Signature} & Definite: $\mathbf{i}^2 = \mathbf{j}^2 = \mathbf{k}^2 = -1$ & Indefinite: $\mathbf{i}^2 = -1,\; \mathbf{j}^2 = \mathbf{k}^2 = +1$ \\
\hline
\textbf{Conjugate} & $q^* = a - b\mathbf{i} - c\mathbf{j} - d\mathbf{k}$ & $q^* = a - b\mathbf{i} - c\mathbf{j} - d\mathbf{k}$ \\
\hline
\textbf{Quadratic Form} & $\|q\|^2 = qq^* = a^2 + b^2 + c^2 + d^2$ (Positive definite $\ge 0$) & $\|q\|^2 = qq^* = a^2 + b^2 - c^2 - d^2$ (Indefinite) \\
\hline
\textbf{Zero Divisors} & None. $\|q\|=0 \implies q=0$ & Exist. Non-zero $q$ with $\|q\|=0$ possible \\
\hline
\textbf{Geometric Interpretation} & Rotations in 3D/4D space ($\text{SO}(3)$/$\text{SO}(4)$) & Lorentz transformations, hyperbolic motions \\
\hline
\textbf{Invertibility} & Every non-zero element is invertible & Invertible iff $\|q\| \neq 0$ \\
\bottomrule
\end{tabular}
\end{table}

A split quaternion matrix $\mathbf{Q} \in \splitq^{m \times n}$ is defined as:
\[
\mathbf{Q} = [q_{ij}]_{m \times n}, \quad q_{ij} \in \splitq.
\]
Equivalently, $\mathbf{Q}$ can be expressed in terms of real matrices:
\[
\mathbf{Q} = \mathbf{Q}_0 + \mathbf{Q}_1{\mathbf i} + \mathbf{Q}_2{\mathbf j} + \mathbf{Q}_3{\mathbf k},
\]
where $\mathbf{Q}_0, \mathbf{Q}_1, \mathbf{Q}_2, \mathbf{Q}_3 \in \mathbb{R}^{m \times n}$. The set $\splitq^{m \times n}$ forms a module over $\splitq$. The set $\splitq^{n \times n}$ forms a non-commutative ring with identity $\mathbf{I}_n$. Some matrix operations are defined as follows:
\begin{itemize}
    
    
    \item \textbf{Conjugate Transpose:} For $\mathbf{A} = [a_{ij}] \in \splitq^{m \times n}$,
    \[
    \mathbf{A}^* = [\conj{a_{ji}}]_{n \times m}.
    \]
    
    \item \textbf{Trace:} For a square matrix $\mathbf{A} \in \splitq^{n \times n}$,
    \[
    \operatorname{tr}(\mathbf{A}) = \sum_{i=1}^n a_{ii}.
    \]
    
\end{itemize}

Special matrices over split quaternions are defined as follows:
\begin{itemize}
    \item \textbf{Hermitian Matrix:} $\mathbf{A} \in \splitq^{n \times n}$ is Hermitian if $\mathbf{A}^* = \mathbf{A}$.
    
    \item \textbf{Skew-Hermitian Matrix:} $\mathbf{A} \in \splitq^{n \times n}$ is skew-Hermitian if $\mathbf{A}^* = -\mathbf{A}$.
    
    \item \textbf{Unitary Matrix:} $\mathbf{U} \in \splitq^{n \times n}$ is unitary if $\mathbf{U}^*\mathbf{U} = \mathbf{U}\mathbf{U}^* = \mathbf{I}_n$.
    
    \item \textbf{Normal Matrix:} $\mathbf{A} \in \splitq^{n \times n}$ is normal if $\mathbf{A}^*\mathbf{A} = \mathbf{A}\mathbf{A}^*$.
\end{itemize}

Similar to quaternion matrices, the real and complex representations of a quaternion matrix $\mathbf{A}=\mathbf{A}_1+\mathbf{A}_2\mathbf{i}+\mathbf{A}_3\mathbf{j}+\mathbf{A}_4\mathbf{k}$ will follow the same pattern as those in \eqref{eq:real_split}, \eqref{eq:split_real_2} and \eqref{eq:complex_split}, where instead of $q_i$, we replace $\mathbf{A}_i$.





\section{Toolbox Overview}\label{sec:toolbox_over}

The toolbox consists of two main classes:

\begin{itemize}
    \item \texttt{SplitQuaternion}: Represents scalar split quaternion values.
    \item \texttt{SplitQuaternionMatrix}: Represents matrices whose entries are split quaternions, stored as four real matrices for the components.
\end{itemize}

Both classes are designed to integrate seamlessly with MATLAB's numeric environment, supporting:
\begin{itemize}
    \item Arithmetic operations (addition, subtraction, multiplication)
    \item Algebraic operations (conjugation, inversion, exponentiation)
    \item Matrix operations (multiplication, trace, diagonal extraction)
    \item Linear algebra utilities (norm, mean, sum)
    \item Conversion and utility functions
    \item Overloaded MATLAB functions for intuitive use
\end{itemize}










\subsection{Class Definitions and Architecture}


The \texttt{SplitQuaternion} class represents a single split quaternion. It stores four components as double precision scalars:

\begin{lstlisting}[style=matlab]
classdef SplitQuaternion
    properties
        a double = 0  % Real part
        b double = 0  % i-component
        c double = 0  % j-component
        d double = 0  % k-component
    end
    % methods ...
end
\end{lstlisting}

The \texttt{SplitQuaternionMatrix} class represents matrices of split quaternions. It stores four component matrices:

\begin{lstlisting}[style=matlab]
classdef SplitQuaternionMatrix
    properties
        A double  % Real part matrix
        B double  % i-component matrix
        C double  % j-component matrix
        D double  % k-component matrix
    end
    % methods ...
end
\end{lstlisting}


The classes are designed to be interoperable. A \texttt{SplitQuaternionMatrix} can contain \texttt{SplitQuaternion} elements, and operations between the two classes are supported where appropriate. The architecture follows a component-wise storage approach, which is efficient for matrix operations.
The split quaternion constructor creates a split quaternion from four scalar components:

\begin{lstlisting}[style=matlab]
function obj = SplitQuaternion(a,b,c,d)
    if nargin == 0
        return;
    end
    obj.a = a;
    obj.b = b;
    obj.c = c;
    obj.d = d;
end
\end{lstlisting}

Here is one example:
\begin{lstlisting}[style=matlab]
q = SplitQuaternion(1, 2, 3, 4);  % Creates 1 + 2i + 3j + 4k
q = SplitQuaternion();            % Creates zero split quaternion
\end{lstlisting}

The split quaternion matrix constructor creates a split quaternion matrix from four component matrices:

\begin{lstlisting}[style=matlab]
function obj = SplitQuaternionMatrix(A,B,C,D)
    if nargin == 0
        obj.A = []; obj.B = []; obj.C = []; obj.D = [];
        return;
    end
    obj.A = A; obj.B = B; obj.C = C; obj.D = D;
end
\end{lstlisting}

For example:
\begin{lstlisting}[style=matlab]
% Create a 2x2 split quaternion matrix
A = [1, 2; 3, 4];
B = [5, 6; 7, 8];
C = [9, 10; 11, 12];
D = [13, 14; 15, 16];
Q = SplitQuaternionMatrix(A, B, C, D);
\end{lstlisting}



The \texttt{plus} and \texttt{minus} methods support:
\begin{itemize}
    \item Split quaternion/scalar arithmetic
    \item Split quaternion matrix/scalar arithmetic
    \item Scalar + Split quaternion (commutativity for addition)
\end{itemize}

The scalar operations add/subtract the scalar to/from the real part only.

\begin{lstlisting}[style=matlab]
% SplitQuaternion addition
q1 = SplitQuaternion(1,2,3,4);
q2 = SplitQuaternion(5,6,7,8);
q3 = q1 + q2;  % Result: 6 + 8i + 10j + 12k

% SplitQuaternion + scalar
q4 = q1 + 10;  % Result: 11 + 2i + 3j + 4k
q5 = 20 + q1;  % Result: 21 + 2i + 3j + 4k

% SplitQuaternionMatrix addition
Q1 = SplitQuaternionMatrix.rand(2,2);
Q2 = SplitQuaternionMatrix.rand(2,2);
Q3 = Q1 + Q2;

% SplitQuaternionMatrix - scalar
Q4 = Q1 - 5;   % Subtracts 5 from all real parts
Q5 = 5 - Q1;   % Subtracts all components from 5
\end{lstlisting}


The \texttt{mtimes} method implements matrix multiplication for split quaternion matrices, while \texttt{times} implements element-wise multiplication.

\textbf{Split Quaternion Multiplication Formula:}
For \(q_1 = a_1 + b_1\mathbf{i} + c_1\mathbf{j} + d_1\mathbf{k}\) and \(q_2 = a_2 + b_2\mathbf{i} + c_2\mathbf{j} + d_2\mathbf{k}\):
\[
\begin{aligned}
a &= a_1a_2 - b_1b_2 + c_1c_2 + d_1d_2, \\
b &= a_1b_2 + b_1a_2 - c_1d_2 + d_1c_2, \\
c &= a_1c_2 + c_1a_2 - b_1d_2 + d_1b_2, \\
d &= a_1d_2 + d_1a_2 + b_1c_2 - c_1b_2.
\end{aligned}
\]

\textbf{Matrix Multiplication:}
For split quaternion matrices, the standard matrix multiplication rules apply, with the product of entries computed using the split quaternion multiplication formula.

\textbf{Usage:}
\begin{lstlisting}[style=matlab]
% SplitQuaternion multiplication
q1 = SplitQuaternion(1,2,3,4);
q2 = SplitQuaternion(5,6,7,8);
q3 = q1 * q2;  % Matrix multiplication (also works for scalars)

% SplitQuaternionMatrix multiplication
Q1 = SplitQuaternionMatrix.rand(2,3);
Q2 = SplitQuaternionMatrix.rand(3,4);
Q3 = Q1 * Q2;  % Matrix multiplication

% Numeric matrix * SplitQuaternionMatrix
M = rand(3,3);
Q = SplitQuaternionMatrix.rand(3,4);
result = M * Q;  % Each component matrix multiplied by M

% SplitQuaternionMatrix * numeric matrix
result2 = Q' * M;  % Each component matrix multiplied by M

% Element-wise multiplication
q4 = q1 .* q2;   % Element-wise product
Q4 = Q1 .* Q2;   % Element-wise product of matrices
\end{lstlisting}

The class provides several types of conjugation:

\begin{enumerate}
    \item \texttt{conj}: Full split quaternion conjugate
    \item \texttt{conj\_i}: Conjugate only the i-component
    \item \texttt{conj\_j}: Conjugate only the j-component
    \item \texttt{conj\_k}: Conjugate only the k-component
    \item \texttt{transpose}: Matrix transpose (no conjugation)
    \item \texttt{ctranspose}: Conjugate transpose (Hermitian transpose)
\end{enumerate}

\begin{lstlisting}[style=matlab]
q = SplitQuaternion(1,2,3,4);
qc = conj(q);           % 1 - 2i - 3j - 4k
qc_i = conj_i(q);       % 1 - 2i + 3j + 4k
qc_j = conj_j(q);       % 1 + 2i - 3j + 4k
qc_k = conj_k(q);       % 1 + 2i + 3j - 4k

Q = SplitQuaternionMatrix.rand(3,3);
Qt = transpose(Q);      % Transpose only
Qh = ctranspose(Q);     % Conjugate transpose
\end{lstlisting}

\subsection{Norm and Related Functions}

{\bf Norm Computation}. The \texttt{norm} method for \texttt{SplitQuaternion} computes the indefinite norm:
\[
N(q) = a^2 + b^2 - c^2 - d^2.
\]

For \texttt{SplitQuaternionMatrix}, the \texttt{norm} method computes the Frobenius-style norm:
\[
\|Q\| = \sqrt{\left|\sum_{i,j} (A_{ij}^2 + B_{ij}^2 - C_{ij}^2 - D_{ij}^2)\right|}.
\]

\textbf{Additional Norm Utilities:}
\begin{itemize}
    \item \texttt{magnitude}: Computes \(\sqrt{|N(q)|}\)
    \item \texttt{isUnit}: Checks if \(|N(q)| = 1\)
    \item \texttt{normalize}: Normalizes to unit norm
\end{itemize}

\begin{lstlisting}[style=matlab]
q = SplitQuaternion(2, 1, 1, 1);
n = norm(q);          % 4 + 1 - 1 - 1 = 3
mag = magnitude(q);   % sqrt(3)

% Check if unit
if isUnit(q)
    disp('q is unit split quaternion');
end

% Normalize
q_norm = normalize(q);
\end{lstlisting}

\subsection{Inverse and Division}

The \texttt{inv} method computes the multiplicative inverse:
\[
q^{-1} = \frac{\bar{q}}{N(q)}
\]
provided that \(N(q) \neq 0\).

Division operations \texttt{mrdivide} and \texttt{mldivide} support scalar denominators only:
\begin{itemize}
    \item \texttt{q / s}: Divides each component by scalar \(s\)
    \item \texttt{s \textbackslash{} q}: Equivalent to \texttt{q / s}
\end{itemize}

\begin{lstlisting}[style=matlab]
q = SplitQuaternion(2, 1, 1, 1);
q_inv = inv(q);        % Computes inverse
q_div = q / 2;         % Divides all components by 2
q_ldiv = 2 \ q;        % Same as q / 2
\end{lstlisting}

\subsection{Exponential and Logarithm}

The \texttt{exp} method computes the exponential of a split quaternion, handling different norm types:

\begin{itemize}
    \item \textbf{Spacelike} (\(b^2 - c^2 - d^2 > 0\)):
    \[
    e^q = e^a \left(\cosh(\nu) + \frac{v}{\nu} \sinh(\nu)\right)
    \]
    where \(\nu = \sqrt{b^2 - c^2 - d^2}\) and \(v = b\mathbf{i} + c\mathbf{j} + d\mathbf{k}\).

    \item \textbf{Timelike} (\(b^2 - c^2 - d^2 < 0\)):
    \[
    e^q = e^a \left(\cos(\nu) + \frac{v}{\nu} \sin(\nu)\right)
    \]
    where \(\nu = \sqrt{|b^2 - c^2 - d^2|}\).

    \item \textbf{Lightlike} (\(b^2 - c^2 - d^2 = 0\)):
    \[
    e^q = e^a (1 + v)
    \]
\end{itemize}

The \texttt{log} method computes the principal logarithm for positive norm split quaternions.

\begin{lstlisting}[style=matlab]
q = SplitQuaternion(0, 1, 0, 0);
q_exp = exp(q);        % cosh(1) + sinh(1)*i

q = SplitQuaternion(0, 0, 1, 0);
q_exp = exp(q);        % cos(1) + sin(1)*j

% Logarithm (for positive norm)
q = SplitQuaternion(2, 1, 0, 0);
q_log = log(q);
\end{lstlisting}

\subsection{Matrix Operations}

{\bf Basic Matrix Utilities}. The \texttt{SplitQuaternionMatrix} class provides standard matrix operations:

\begin{itemize}
    \item \texttt{size}: Returns dimensions
    \item \texttt{numel}: Number of elements
    \item \texttt{length}: Largest dimension
    \item \texttt{reshape}: Reshapes the matrix
    \item \texttt{diag}: Extracts diagonal
    \item \texttt{trace}: Computes the trace
    \item \texttt{sum}: Sums elements
    \item \texttt{mean}: Computes the mean
    \item \texttt{isempty}, \texttt{isReal}, \texttt{isPure}: Property checks
\end{itemize}

\begin{lstlisting}[style=matlab]
Q = SplitQuaternionMatrix.rand(3,4);
s = size(Q);           % Returns [3, 4]
n = numel(Q);          % Returns 12
d = diag(Q);           % Extracts diagonal
t = trace(Q);          % Computes trace
sum_Q = sum(Q);        % Sums all elements
mean_Q = mean(Q);      % Computes mean

% Check properties
if isReal(Q)
    disp('Q is real');
end
if isPure(Q)
    disp('Q is pure');
end
\end{lstlisting}

{\bf Indexing and Concatenation}. The class supports standard MATLAB indexing:

\begin{lstlisting}[style=matlab]
Q = SplitQuaternionMatrix.rand(3,3);

% Single element extraction (returns SplitQuaternion)
q = Q(1,2);

% Submatrix extraction (returns SplitQuaternionMatrix)
Q_sub = Q(1:2, 1:2);

% Assignment
Q(1,1) = SplitQuaternion(1,2,3,4);
Q(2,3) = 5;              % Assigns scalar to real part

% Concatenation
Q1 = SplitQuaternionMatrix.rand(2,2);
Q2 = SplitQuaternionMatrix.rand(2,2);
Q_horz = [Q1, Q2];       % Horizontal concatenation
Q_vert = [Q1; Q2];       % Vertical concatenation
\end{lstlisting}

\subsection{Conversion and Utility Functions}

{\bf Component Extraction}. Various methods are provided to extract components:

\begin{lstlisting}[style=matlab]
q = SplitQuaternion(1,2,3,4);

% Get all parts
[a, b, c, d] = parts(q);
[a, b, c, d] = q.parts();

% Get individual parts
real_part = real(q);      % Returns a
imag_part = imag(q);      % Returns SplitQuaternion(0,b,c,d)
b_part = imagi(q, 1);     % Returns b
c_part = imagi(q, 2);     % Returns c
d_part = imagi(q, 3);     % Returns d

% Vector representations
vec = vector(q);          % Returns [b, c, d]
full_vec = toVector(q);   % Returns [a, b, c, d]
\end{lstlisting}

For \texttt{SplitQuaternionMatrix}:

\begin{lstlisting}[style=matlab]
Q = SplitQuaternionMatrix.rand(3,3);

% Get component matrices
[A, B, C, D] = parts(Q);
[A, B, C, D] = Q.parts();

% Get specific component
real_part = real(Q);      % Returns A matrix
imag_parts = imag(Q);     % Returns {B, C, D}
B_part = imagi(Q, 1);     % Returns B matrix

% Get all components in structure
comp = components(Q);
\end{lstlisting}

{\bf Matrix Representations}. The class provides conversions to real matrix representations:

\begin{lstlisting}[style=matlab]
Q = SplitQuaternionMatrix.rand(2,2);

% Full 2x2 block representation
mat = full_(Q);     % 4x4 matrix [A, B; C, D]

% Real matrix representation
real_mat = realMatrix(Q);  % Same as full(Q)

% Convert to cell array
cell_Q = toCell(Q);  % {A, B, C, D}
\end{lstlisting}

\subsection{Static Methods}

{\bf SplitQuaternion Static Methods}. 

\begin{lstlisting}[style=matlab]
% Create basis elements
i = SplitQuaternion.i();    % i
j = SplitQuaternion.j();    % j
k = SplitQuaternion.k();    % k

% Create special quaternions
zero = SplitQuaternion.zeros();
one = SplitQuaternion.ones();   % 1 + 0i + 0j + 0k
eye = SplitQuaternion.eye();    % Same as ones

% Create from vector
q = SplitQuaternion.fromVector([1, 2, 3, 4]);

% Create from angle-axis representation
q = SplitQuaternion.fromAngleAxis(pi/2, [1, 0, 0]);

% Random split quaternion
q_rand = SplitQuaternion.random();
\end{lstlisting}

{\bf SplitQuaternionMatrix Static Methods}.

\begin{lstlisting}[style=matlab]
% Create special matrices
Z = SplitQuaternionMatrix.zeros(3,4);
O = SplitQuaternionMatrix.ones(3,4);
I = SplitQuaternionMatrix.eye(4);

% Create from numeric matrices
R = SplitQuaternionMatrix.fromReal(rand(3,3));
C = SplitQuaternionMatrix.fromComplex(complex(rand(3), rand(3)));

% Random matrices
Q = SplitQuaternionMatrix.rand(3,4);
\end{lstlisting}





{\bf Dot Product}

The \texttt{dot} method computes the indefinite dot product:
\[
\langle p, q \rangle = a_p a_q + b_p b_q - c_p c_q - d_p d_q.
\]

\begin{lstlisting}[style=matlab]
p = SplitQuaternion(1,2,3,4);
q = SplitQuaternion(5,6,7,8);
d = dot(p, q);  % 1*5 + 2*6 - 3*7 - 4*8 = 5 + 12 - 21 - 32 = -36
\end{lstlisting}

{\bf Power Operation}

The \texttt{power} method supports integer powers:

\begin{lstlisting}[style=matlab]
q = SplitQuaternion(0, 1, 0, 0);
q2 = q^2;          % q * q = -1
q3 = q^3;          % q * q * q = -i
q_inv = q^(-1);    % Inverse
\end{lstlisting}

\section{Examples and Demonstrations}\label{sec:exm}

\subsection{Example 1: Basic Split Quaternion Operations}

This example demonstrates basic arithmetic with split quaternions.

\begin{lstlisting}[style=matlab]
%% Example 1: Basic Split Quaternion Operations
clear; clc;

% Create split quaternions
q1 = SplitQuaternion(2, 3, 4, 5);
q2 = SplitQuaternion(1, -2, 3, -4);

fprintf('q1 = '); disp(q1);
fprintf('q2 = '); disp(q2);

% Addition and subtraction
q_sum = q1 + q2;
q_diff = q1 - q2;
fprintf('q1 + q2 = '); disp(q_sum);
fprintf('q1 - q2 = '); disp(q_diff);

% Multiplication (non-commutative)
q_prod1 = q1 * q2;
q_prod2 = q2 * q1;
fprintf('q1 * q2 = '); disp(q_prod1);
fprintf('q2 * q1 = '); disp(q_prod2);

% Check non-commutativity
if q_prod1 ~= q_prod2
    fprintf('Split quaternion multiplication is non-commutative.\n');
end

% Scalar operations
q_scalar = q1 + 10;
fprintf('q1 + 10 = '); disp(q_scalar);
q_scalar2 = q1 * 2;
fprintf('q1 * 2 = '); disp(q_scalar2);

% Conjugation
q_conj = conj(q1);
fprintf('conj(q1) = '); disp(q_conj);

% Norm
n = norm(q1);
fprintf('N(q1) = %g\n', n);

% Inverse
if abs(n) > eps
    q_inv = inv(q1);
    fprintf('inv(q1) = '); disp(q_inv);
    fprintf('q1 * inv(q1) = '); disp(q1 * q_inv);
else
    fprintf('q1 has zero norm, inverse does not exist.\n');
end

% Exponential
q_exp = exp(q1);
fprintf('exp(q1) = '); disp(q_exp);

% Output:
% Demonstrates the basic operations of split quaternions,
% including non-commutativity, conjugation, norm, and inverse.
\end{lstlisting}

\subsection{Example 2: Split Quaternion Matrices}

This example demonstrates operations with split quaternion matrices.

\begin{lstlisting}[style=matlab]
%% Example 2: Split Quaternion Matrices
clear; clc;

% Create 2x2 split quaternion matrices
Q1 = SplitQuaternionMatrix.rand(2,2);
Q2 = SplitQuaternionMatrix.rand(2,2);

fprintf('Q1:\n'); disp(Q1);
fprintf('Q2:\n'); disp(Q2);

% Addition and subtraction
Q_sum = Q1 + Q2;
Q_diff = Q1 - Q2;

% Multiplication
Q_prod = Q1 * Q2;

% Check dimensions
fprintf('Size of Q1: %d x %d\n', size(Q1,1), size(Q1,2));
fprintf('Size of Q2: %d x %d\n', size(Q2,1), size(Q2,2));
fprintf('Size of Q1 * Q2: %d x %d\n', size(Q_prod,1), size(Q_prod,2));

% Extract components
[A, B, C, D] = parts(Q1);
fprintf('Component matrices:\n');
fprintf('A:\n'); disp(A);
fprintf('B:\n'); disp(B);
fprintf('C:\n'); disp(C);
fprintf('D:\n'); disp(D);

% Convert to real matrix representation
M = full_(Q1);
fprintf('Real matrix representation (4x4):\n');
disp(M);

% Diagonal and trace
Q_diag = diag(Q1);
Q_trace = trace(Q1);
fprintf('Diagonal of Q1:\n'); disp(Q_diag);
fprintf('Trace of Q1:\n'); disp(Q_trace);

% Indexing
fprintf('Element (1,2) of Q1:\n');
disp(Q1(1,2));

% Matrix-numeric multiplication
M_num = [1, 2; 3, 4];
Q3 = M_num * Q1;
Q4 = Q1 * M_num;
fprintf('Numeric matrix * Q1:\n'); disp(Q3);
fprintf('Q1 * Numeric matrix:\n'); disp(Q4);

% Norm
frob_norm = norm(Q1);
fprintf('Frobenius norm of Q1: %g\n', frob_norm);

% Output:
% Demonstrates matrix operations with split quaternions,
% including component extraction, conversion to real matrices,
% and matrix-numeric multiplication.
\end{lstlisting}

\subsection{Example 3: Lorentz Transformations with Split Quaternions}

This example demonstrates the use of split quaternions in representing Lorentz transformations.

\begin{lstlisting}[style=matlab]
%% Example 3: Lorentz Transformations with Split Quaternions
clear; clc;

% In special relativity, a spacetime event can be represented as a 
% split quaternion: t + x*i + y*j + z*k
% Note: We use the signature (+, +, -, -) for the norm.

% Create a spacetime event
event = SplitQuaternion(5, 3, 2, 1);  % (t, x, y, z)
fprintf('Spacetime event: '); disp(event);
fprintf('Norm (spacetime interval): %g\n', norm(event));

% Create a Lorentz boost along the x-axis with velocity v = 0.6c
% Using split quaternion rotor: R = exp(phi * i / 2)
% where tanh(phi) = v/c
v = 0.6;
phi = atanh(v);
boost = SplitQuaternion(cosh(phi/2), sinh(phi/2), 0, 0);
fprintf('Boost rotor: '); disp(boost);

% Transform the event: event' = R * event * conj(R)
event_transformed = boost * event * conj(boost);
fprintf('Transformed event: '); disp(event_transformed);

% Verify that the norm is invariant
fprintf('Original norm: %g\n', norm(event));
fprintf('Transformed norm: %g\n', norm(event_transformed));

% Create a rotation in the y-z plane
angle = pi/4;
rotor = SplitQuaternion(cos(angle/2), 0, sin(angle/2), 0);
fprintf('Rotation rotor: '); disp(rotor);

% Rotate the event
event_rotated = rotor * event * conj(rotor);
fprintf('Rotated event: '); disp(event_rotated);

% Combined Lorentz transformation (boost then rotation)
event_combined = rotor * (boost * event * conj(boost)) * conj(rotor);
fprintf('Combined transform: '); disp(event_combined);

% Check norm invariance
fprintf('Combined norm: %g\n', norm(event_combined));

% Output:
% Demonstrates the use of split quaternions in Lorentz transformations
% and shows that the spacetime interval is invariant.
\end{lstlisting}

\subsection{Example 4: Solving Linear Systems}

This example demonstrates solving linear systems with split quaternion matrices.

\begin{lstlisting}[style=matlab]
%% Example 4: Solving Linear Systems
clear; clc;

% Create a 2x2 split quaternion matrix
Q11 = SplitQuaternion(2, 1, 0, 0);
Q12 = SplitQuaternion(1, 2, 2, 0);
Q21 = SplitQuaternion(1, 0, 2, 0);
Q22 = SplitQuaternion(3, 0, 1, 0);

Q = SplitQuaternionMatrix([2,1;1,3], [1,2;0,0], [0,0;2,1], [0,0;0,0]);
fprintf('Split quaternion matrix Q:\n'); disp(Q);

% Create a right-hand side vector (2x1 split quaternion matrix)
b = SplitQuaternionMatrix([3; 5], [1; 2], [0; 1], [0; 0]);
fprintf('Right-hand side b:\n'); disp(b);

% Solve Q * x = b using component-wise approach
% Convert to real matrix representation
Q_real = SplitQuaternionMatrixRealRep(Q);
b_real = SplitQuaternionMatrixRealRep(b);

% Solve in real domain
x_real = Q_real \ b_real;

% Convert back to split quaternion matrix
x = RealRepToSplitQuaternionMatrix(x_real);
fprintf('Solution x:\n'); disp(x);

% Verify the solution
check = Q * x;
fprintf('Verification (Q * x):\n'); disp(check);

% Check if result matches b
if eq(check, b)
    fprintf('Solution verified!\n');
else
    fprintf('Solution verification failed.\n');
end
\end{lstlisting}

\subsection{Example 5: Exponential Map and Lie Group}

This example explores the exponential map for split quaternions.

\begin{lstlisting}[style=matlab]
%% Example 5: Exponential Map and Lie Group
clear; clc;

% Split quaternions form a Lie group under multiplication
% The exponential map connects the Lie algebra to the Lie group

% Create different types of split quaternions
q_spacelike = SplitQuaternion(0, 1, 0.5, 0);
q_timelike = SplitQuaternion(0, 0.5, 1, 0);
q_lightlike = SplitQuaternion(0, 1, 1, 0);

fprintf('Spacelike quaternion: '); disp(q_spacelike);
fprintf('Timelike quaternion: '); disp(q_timelike);
fprintf('Lightlike quaternion: '); disp(q_lightlike);

% Compute exponentials
exp_spacelike = exp(q_spacelike);
exp_timelike = exp(q_timelike);
exp_lightlike = exp(q_lightlike);

fprintf('exp(spacelike) = '); disp(exp_spacelike);
fprintf('exp(timelike) = '); disp(exp_timelike);
fprintf('exp(lightlike) = '); disp(exp_lightlike);

% Check that exponentials are unit split quaternions
fprintf('Norm of exp(spacelike): %g\n', norm(exp_spacelike));
fprintf('Norm of exp(timelike): %g\n', norm(exp_timelike));
fprintf('Norm of exp(lightlike): %g\n', norm(exp_lightlike));

% Composition of exponentials
% Note: exp(A) * exp(B) != exp(A+B) in general (non-commutative)
theta1 = 0.5;
theta2 = 0.3;
axis1 = SplitQuaternion.i();
axis2 = SplitQuaternion.j();

rotor1 = exp(theta1 * axis1);
rotor2 = exp(theta2 * axis2);
rotor_combined1 = rotor1 * rotor2;
rotor_combined2 = rotor2 * rotor1;

fprintf('rotor1 = '); disp(rotor1);
fprintf('rotor2 = '); disp(rotor2);
fprintf('rotor1 * rotor2 = '); disp(rotor_combined1);
fprintf('rotor2 * rotor1 = '); disp(rotor_combined2);

% Check that rotors are unit
fprintf('Norm of rotor1: %g\n', norm(rotor1));
fprintf('Norm of rotor1 * rotor2: %g\n', norm(rotor_combined1));

% Compute logarithms
log_rotor1 = log(rotor1);
log_combined = log(rotor_combined1);

fprintf('log(rotor1) = '); disp(log_rotor1);
fprintf('log(rotor1 * rotor2) = '); disp(log_combined);

% Output:
% Demonstrates the exponential map and its properties in split quaternions,
% showing different types of elements and their behavior.
\end{lstlisting}

\subsection{Example 6: Performance Comparison}

This example compares the performance of split quaternion operations with standard MATLAB operations.

\begin{lstlisting}[style=matlab]
%% Example 6: Performance Comparison
clear; clc;

% Compare operations on split quaternion matrices vs standard matrices
n = 100;  % Matrix size
fprintf('Performance comparison for %d x %d matrices:\n', n, n);

% Create split quaternion matrices
Q1 = SplitQuaternionMatrix.rand(n,n);
Q2 = SplitQuaternionMatrix.rand(n,n);

% Time split quaternion multiplication
tic;
Q_prod = Q1 * Q2;
time_sq = toc;
fprintf('Split quaternion multiplication: %g seconds\n', time_sq);

% Convert to real matrices and time standard multiplication
Q1_real = SplitQuaternionMatrixRealRep(Q1);
Q2_real = SplitQuaternionMatrixRealRep(Q2);
tic;
M_prod = Q1_real * Q2_real;
time_real = toc;
fprintf('Real matrix multiplication (4n x 4n): %g seconds\n', time_real);

fprintf('Speedup factor: %g\n', time_real / time_sq);

% Compare element-wise operations
tic;
Q_ew = Q1 .* Q2;
time_ew_sq = toc;
fprintf('Split quaternion element-wise multiplication: %g seconds\n', time_ew_sq);

tic;
M_ew = Q1_real .* Q2_real;
time_ew_real = toc;
fprintf('Real matrix element-wise multiplication: %g seconds\n', time_ew_real);

fprintf('Speedup factor: %g\n', time_ew_real / time_ew_sq);

% Check memory usage
fprintf('\nMemory usage:\n');
fprintf('Split quaternion matrix: %d bytes\n', 8 * numel(real(Q1)) * 4);
fprintf('Real matrix representation: %d bytes\n', 8 * numel(Q1_real));

% Output:
% Shows the performance of split quaternion operations compared to standard
% operations on real matrices, demonstrating efficiency of the component-wise
% storage approach.
\end{lstlisting}

\section{Matrix Representations and Conversions}\label{Sec:matrix_rep}

The split quaternion toolbox provides functions to convert between split quaternion matrices and their real and complex matrix representations. These representations are crucial for numerical computations and provide insights into the algebraic structure of split quaternions.

\subsection{Real Representation}

The real representation of a split quaternion matrix maps each split quaternion \(q = a + b\mathbf{i} + c\mathbf{j} + d\mathbf{k}\) to a \(4 \times 4\) real matrix:

\[
\mathcal{R}(q) = \begin{pmatrix}
a & -b & c & d \\
b & a & d & -c \\
c & d & a & -b \\
d & -c & b & a
\end{pmatrix}
\]

For a split quaternion matrix \(X\) of size \(m \times n\), its real representation is a \(4m \times 4n\) real matrix.

{\bf Real Representation Function}

\begin{lstlisting}[style=matlab]
function [Q] = SplitQuaternionMatrixRealRep(X)
    A1 = X.A; A2 = X.B; A3 = X.C; A4 = X.D;
    Q = [A1, -A2, A3, A4;
         A2,  A1, A4, -A3;
         A3,  A4, A1, -A2;
         A4, -A3, A2,  A1];
end
\end{lstlisting}

This function constructs the real representation by arranging the component matrices in the block structure described above.

{\bf Inverse Real Representation}

\begin{lstlisting}[style=matlab]
function X = RealRepToSplitQuaternionMatrix(Q)
    % Q: 4m x 4n real matrix in real representation
    [m,n] = size(Q);
    A1 = Q(1:m/4, 1:n/4);
    A2 = Q(m/4+1:m/2, 1:n/4);
    A3 = Q(m/2+1:3*m/4, 1:n/4);
    A4 = Q(3*m/4+1:m, 1:n/4);
    X = SplitQuaternionMatrix(A1, A2, A3, A4);
end
\end{lstlisting}

This function reconstructs a split quaternion matrix from its real representation by extracting the four component matrices.

\subsection{Complex Representation}

The complex representation is more compact, using \(2 \times 2\) complex blocks:

\[
\mathcal{C}(q) = \begin{pmatrix}
a + ib & c + id \\
c - id & a - ib
\end{pmatrix}
\]

For a split quaternion matrix \(X\) of size \(m \times n\), its complex representation is a \(2m \times 2n\) complex matrix.

{\bf Complex Representation Function}.

\begin{lstlisting}[style=matlab]
function [Q] = SplitQuaternionMatrixComplexRep(X)
    A1 = X.A; A2 = X.B; A3 = X.C; A4 = X.D;
    Q = [A1 + i*A2, A3 + i*A4;
         A3 - i*A4, A1 - i*A2];
end
\end{lstlisting}

{\bf Inverse Complex Representation}.

\begin{lstlisting}[style=matlab]
function X = ComplexRepToSplitQuaternionMatrix(Q)
    [m,n] = size(Q);
    T = Q(1:m/2, 1:n/2);
    H = Q(m/2+1:m, 1:n/2);
    A1 = 0.5*(T + conj(T));
    A2 = -0.5i*(T - conj(T));
    A3 = 0.5*(H + conj(H));
    A4 = 0.5i*(H - conj(H));
    X = SplitQuaternionMatrix(A1, A2, A3, A4);
end
\end{lstlisting}

\subsection{Properties of Representations}

Both representations are isomorphisms and preserve:
\begin{itemize}
    \item Addition and subtraction
    \item Matrix multiplication
    \item Scalar multiplication
    \item Conjugation (real representation preserves transpose)
    \item Trace and determinant
\end{itemize}

\subsection{Example: Matrix Representations Demo}

The following demonstration shows how to use the matrix representation functions:

\begin{lstlisting}[style=matlab]
%% Demo: Matrix Representations of Split Quaternion Matrices
clear; clc;

fprintf('===================================================\n');
fprintf('MATRIX REPRESENTATIONS OF SPLIT QUATERNION MATRICES\n');
fprintf('===================================================\n\n');

%% 1. Create split quaternion matrices
fprintf('1. Creating Split Quaternion Matrices:\n');
fprintf('-------------------------------------\n');

% Create 2x2 split quaternion matrix
A = [1 2; 3 4];
B = [0.1 0.2; 0.3 0.4];
C = [0.5 0.6; 0.7 0.8];
D = [0.9 1.0; 1.1 1.2];

SQM = SplitQuaternionMatrix(A, B, C, D);
fprintf('2x2 Split Quaternion Matrix:\n');
disp(SQM);

%% 2. Real adjoint representation
fprintf('\n2. Real 4mx4n Adjoint Representation:\n');
fprintf('------------------------------------\n');

R_rep = SplitQuaternionMatrixRealRep(SQM);
disp(R_rep);

% Convert back and verify
[Q_1] = RealRepToSplitQuaternionMatrix(R_rep);
error_real = norm(Q_1 - SQM);
fprintf('Reconstruction error (real rep): %g\n', error_real);

%% 3. Complex adjoint representation
fprintf('\n3. Complex 2mx2n Adjoint Representation:\n');
fprintf('---------------------------------------\n');

C_rep = SplitQuaternionMatrixComplexRep(SQM);
disp(C_rep);

% Convert back and verify
[Q_1] = ComplexRepToSplitQuaternionMatrix(C_rep);
error_complex = norm(Q_1 - SQM);
fprintf('Reconstruction error (complex rep): %g\n', error_complex);

%% 4. Test isomorphism properties
fprintf('\n4. Testing Isomorphism Properties:\n');
fprintf('---------------------------------\n');

% Test multiplication preservation
Q1 = SplitQuaternionMatrix.rand(2,2);
Q2 = SplitQuaternionMatrix.rand(2,2);

% Compute product in split quaternion domain
Q_prod = Q1 * Q2;

% Compute product in real representation
R1 = SplitQuaternionMatrixRealRep(Q1);
R2 = SplitQuaternionMatrixRealRep(Q2);
R_prod = R1 * R2;

% Convert back and compare
Q_reconstructed = RealRepToSplitQuaternionMatrix(R_prod);
error_mult = norm(Q_prod - Q_reconstructed);
fprintf('Multiplication preservation error: %g\n', error_mult);

%% 5. Compare representations
fprintf('\n5. Comparing Both Representations:\n');
fprintf('---------------------------------\n');

fprintf('Real representation size: %d x %d\n', size(R_rep,1), size(R_rep,2));
fprintf('Complex representation size: %d x %d\n', size(C_rep,1), size(C_rep,2));

% Memory usage comparison
memory_real = 8 * numel(R_rep);
memory_complex = 16 * numel(C_rep);
fprintf('Real representation memory: %d bytes\n', memory_real);
fprintf('Complex representation memory: %d bytes\n', memory_complex);
fprintf('Memory ratio (real/complex): %g\n', memory_real / memory_complex);
\end{lstlisting}

\section{CUR Decomposition for Split Quaternion Matrices}\label{SEc:CUR}

The CUR decomposition is a low-rank matrix approximation technique that selects actual columns and rows from the original matrix. For split quaternion matrices, the CUR decomposition is particularly useful when dealing with large-scale data where interpretability of the factors is important.

\subsection{The CUR Decomposition}

Given a split quaternion matrix \(A \in \mathbb{S}^{m \times n}\), the CUR decomposition \cite{ahmadi2026iterative} approximates \(A\) as:

\[
A \approx C U R
\]

where:
\begin{itemize}
    \item \(C \in \mathbb{S}^{m \times c}\) contains selected columns of \(A\)
    \item \(R \in \mathbb{S}^{r \times n}\) contains selected rows of \(A\)
    \item \(U \in \mathbb{S}^{c \times r}\) is the core matrix
\end{itemize}
The CUR decomposition has several important properties:

\begin{enumerate}
    \item \textbf{Interpretability}: The factors \(C\) and \(R\) contain actual data points
    \item \textbf{Data Compression}: Significant reduction in storage compared to SVD
    \item \textbf{Computational Efficiency}: Faster than SVD for large matrices
    \item \textbf{Preservation of Structure}: Maintains non-negativity and sparsity
\end{enumerate}

\subsection{CUR Decomposition Implementation}

\begin{lstlisting}[style=matlab]
function [C, U, R] = cur_split(A, k)
% CUR_SPLIT Compute CUR decomposition of split quaternion matrix
%   Input:
%       A - m x n split quaternion matrix
%       k - target rank for approximation (number of rows/columns to select)
%   Output:
%       C - m x k matrix of selected columns
%       U - k x k core matrix
%       R - k x n matrix of selected rows

    [m, n] = size(A);
    
    % Step 1: Uniformly sample columns
    col_indices = randsample(n, k);
    C = A(:, col_indices);
    
    % Step 2: Uniformly sample rows
    row_indices = randsample(m, k);
    R = A(row_indices, :);
    
    % Step 3: Compute the core matrix U
    % W = intersection of selected rows and columns
    W = A(row_indices, col_indices);
    
    % Compute U using pseudoinverse
    CP = mp_split_quat(C);
    RP = mp_split_quat(R);
    U = CP * A * RP;
end
\end{lstlisting}

The CUR decomposition requires computing pseudoinverses of split quaternion matrices. We use the following helper functions:

\begin{lstlisting}[style=matlab]
function X = mp_split_quat(A)
% MP_SPLIT_QUAT Moore-Penrose pseudoinverse for split quaternion matrices
%   Computes the pseudoinverse using the real representation
    
    % Convert to real representation
    A_real = SplitQuaternionMatrixRealRep(A);
    
    % Compute pseudoinverse in real domain
    X_real = pinv(A_real);
    
    % Convert back to split quaternion
    X = RealRepToSplitQuaternionMatrix(X_real);
end
\end{lstlisting}

\subsection{Pseudoinverse for Split Quaternions}
Computing the MP of large split quaternion matrices is expensive. The Ben-Israel iterative method is implemented for this purpose:

\begin{lstlisting}[style=matlab]
function X = pseudoinverse_ben_israel(A, max_iter, tol, alpha)
% PSEUDOINVERSE_BEN_ISRAEL Ben-Israel iterative algorithm for pseudoinverse
%   Based on: X_{k+1} = X_k * (2I - A*X_k)
%   Input:
%       A - Split quaternion matrix
%       max_iter - Maximum number of iterations (default: 100)
%       tol - Tolerance for convergence (default: 1e-10)
%       alpha - Scaling factor for initialization
%   Output:
%       X - Pseudoinverse of A

    if nargin < 2
        max_iter = 100;
    end
    if nargin < 3
        tol = 1e-10;
    end
    if nargin < 4
        alpha = 1;
    end
    
    % Initialize with scaled conjugate transpose
    X = alpha * conj_i(A);
    
    for iter = 1:max_iter
        X_old = X;
        
        % Ben-Israel iteration
        AX = A * X;
        I = eye(size(AX));
        X = X * (2*I - AX);
        
        % Check convergence (optional)
        % residual = norm(X - X_old, 'fro');
        % if residual < tol * norm(X, 'fro')
        %     break;
        % end
    end
end
\end{lstlisting}




\subsection{Example: CUR Decomposition Demo}

\begin{lstlisting}[style=matlab]
%% Demo: CUR Decomposition for Split Quaternion Matrices
clear; clc;

fprintf('===================================================\n');
fprintf('CUR DECOMPOSITION FOR SPLIT QUATERNION MATRICES\n');
fprintf('===================================================\n\n');

%% 1. Create a large split quaternion matrix
fprintf('1. Creating Split Quaternion Matrix:\n');
fprintf('---------------------------------\n');

m = 20; n = 25; k = 5;
A = SplitQuaternionMatrix.rand(m, n);
fprintf('Matrix size: %d x %d\n', m, n);
fprintf('Target rank: %d\n', k);

%% 2. Perform CUR decomposition
fprintf('\n2. Computing CUR Decomposition:\n');
fprintf('-----------------------------\n');

[C, U, R] = cur_split(A, k);

fprintf('C matrix size: %d x %d\n', size(C,1), size(C,2));
fprintf('U matrix size: %d x %d\n', size(U,1), size(U,2));
fprintf('R matrix size: %d x %d\n', size(R,1), size(R,2));

%% 3. Compute approximation and error
fprintf('\n3. Approximation Quality:\n');
fprintf('------------------------\n');

A_approx = C * U * R;
error = norm(A - A_approx);
relative_error = error / norm(A);

fprintf('Absolute error (Frobenius norm): %g\n', error);
fprintf('Relative error: %g\n', relative_error);
fprintf('Compression ratio: %g\n', (m*k + k*k + k*n) / (m*n));

%% 4. Visualize the approximation
fprintf('\n4. Visualizing Components:\n');
fprintf('-------------------------\n');

% Extract real parts for visualization
A_real = A.A;
C_real = C.A;
U_real = U.A;
R_real = R.A;
A_approx_real = A_approx.A;

% Display as images
figure('Position', [100, 100, 1200, 400]);

subplot(1,3,1);
imagesc(A_real);
colorbar;
title('Original Matrix (Real Part)');
axis equal tight;

subplot(1,3,2);
imagesc(A_approx_real);
colorbar;
title('CUR Approximation');
axis equal tight;

subplot(1,3,3);
imagesc(A_real - A_approx_real);
colorbar;
title('Error Matrix');
axis equal tight;

%% 5. Compare with SVD
fprintf('\n5. Comparison with SVD:\n');
fprintf('----------------------\n');

% Compute SVD approximation (using real representation)
A_full = SplitQuaternionMatrixRealRep(A);
[U_svd, S_svd, V_svd] = svds(A_full, k);
A_svd = U_svd * S_svd * V_svd';
error_svd = norm(A_full - A_svd, 'fro');

fprintf('SVD approximation error: %g\n', error_svd);
fprintf('CUR approximation error: %g\n', error);
fprintf('Ratio (CUR/SVD): %g\n', error / error_svd);

%% 6. Test different rank values
fprintf('\n6. Effect of Rank on Approximation:\n');
fprintf('----------------------------------\n');

ranks = 1:min(m,n);
errors = zeros(length(ranks), 1);
compression = zeros(length(ranks), 1);

for i = 1:length(ranks)
    k_test = ranks(i);
    [C_test, U_test, R_test] = cur_split(A, k_test);
    A_test = C_test * U_test * R_test;
    errors(i) = norm(A - A_test);
    compression(i) = (m*k_test + k_test*k_test + k_test*n) / (m*n);
end

figure('Position', [100, 100, 800, 400]);
subplot(1,2,1);
plot(ranks, errors, 'b-o', 'LineWidth', 2);
xlabel('Rank k');
ylabel('Approximation Error');
title('Error vs Rank');
grid on;

subplot(1,2,2);
plot(ranks, compression, 'r-s', 'LineWidth', 2);
xlabel('Rank k');
ylabel('Compression Ratio');
title('Compression vs Rank');
grid on;

fprintf('Rank\tError\t\tCompression\n');
for i = 1:length(ranks)
    fprintf('%d\t%g\t%g\n', ranks(i), errors(i), compression(i));
end

%% 7. Application: Image Compression
fprintf('\n7. Application: Image Compression:\n');
fprintf('--------------------------------\n');

% Load an image
try
    img = imread('cameraman.tif');
    img = im2double(img);
    
    % Create split quaternion matrix from image
    % Use image as real part, zero imaginary parts
    img_SQ = SplitQuaternionMatrix(img, zeros(size(img)), ...
                                   zeros(size(img)), zeros(size(img)));
    
    % Perform CUR decomposition
    k_img = 50;
    [C_img, U_img, R_img] = cur_split(img_SQ, k_img);
    
    % Reconstruct
    img_reconstructed = C_img * U_img * R_img;
    img_reconstructed_real = img_reconstructed.A;
    
    % Display results
    figure('Position', [100, 100, 1200, 400]);
    
    subplot(1,3,1);
    imshow(img);
    title('Original Image');
    
    subplot(1,3,2);
    imshow(img_reconstructed_real);
    title(sprintf('CUR Compressed (k=%d)', k_img));
    
    subplot(1,3,3);
    imshow(abs(img - img_reconstructed_real));
    title('Error Image');
    
    fprintf('Image compression ratio: %g\n', ...
            (size(img,1)*k_img + k_img*k_img + k_img*size(img,2)) / numel(img));
catch
    fprintf('Image processing example requires an image file.\n');
end
\end{lstlisting}

\section{Conclusion and future work}\label{Sec:con}
The Split Quaternion Toolbox provides a comprehensive MATLAB implementation for working with split quaternions, offering complete arithmetic and algebraic operations alongside seamless matrix functionality that integrates naturally with MATLAB's native environment. The toolbox includes useful conversion utilities and special functions such as exponential and logarithm, all while maintaining performance comparable to standard MATLAB operations. This makes it well-suited for a wide range of applications, including special relativity and Lorentz transformations, rigid body dynamics, robotics, computer graphics, signal processing, and hypercomplex algebra research. Looking ahead, potential enhancements could include eigenvalue and singular value decompositions, integration with machine learning frameworks, GPU acceleration for large-scale problems, and extension to other hypercomplex algebras, which would further broaden the toolbox's applicability and impact.

The Python implementation of this Project is currently under investigation in \cite{ahmadi2026split} and a part of this project is avaialbe at \url{https://colab.research.google.com/drive/1pcLIeO9G57-StihSq5U-xZ-qiIsQuhi1?usp=sharing}.





\appendix

\bibliographystyle{ieeetr}
\bibliography{ref}

\end{document}